\documentclass[12pt]{article}
\usepackage{amsmath}
\usepackage[pdftex]{graphicx}

\begin{document}
\begin{center} \Large{{\bf{Self-Exciting Random Evolutions (SEREs) and Their Applications (Version 2)}}\\

Anatoliy Swishchuk\footnote {The author thanks to two unanimous referees for their valuable comments and remarks with respect to the first version of this working paper (Available at SSRN: https://ssrn.com/abstract=5055075). The work on this working paper is still in progress.}\\

Department Mathematics \& Statistics\\
University of Calgary\\
Calgary, Canada T2N 1N4\\
aswish@ucalgary.ca}

\end{center}

\begin{abstract}
This paper is devoted to the study of a new class of random evolutions (RE), so-called self-exciting random evolutions (SEREs), and their applications. We also introduce a new random process $x(t)$ such that it is based on  a superposition of a Markov chain $x_n$ and a Hawkes process $N(t),$ i.e., $x(t):=x_{N(t)}.$ We call this process self-walking imbedded semi-Hawkes process (Swish Process or SwishP). Then the self-exciting REs (SEREs) can be constructed in similar way as, e.g., semi-Markov REs, but instead of semi-Markov process $x(t)$ we have SwishP.  We give classifications and examples of self-exciting REs (SEREs). Then we consider two limit theorems for SEREs such as  averaging (Theorem 1) and diffusion approximation (Theorem 2). Applications of SEREs are devoted to the so-called self-exciting traffic/transport process and self-exciting summation on a Markov chain, which are examples of continuous and discrete SERE. From these processes we can construct many other self-exciting processes, e.g., such as impulse traffic/transport process, self-exciting risk process, general compound Hawkes process for a stock price, etc. We present averaged and diffusion approximation of self-exciting processes. The novelty of the paper associated with new models, such as $x(t)$ and SERE, and also new features of SEREs and their many applications, namely, self-exciting and clustering effects.

\end{abstract}


\section{Introduction}

The purpose of the present paper is to introduce and to study a new class of Random Evolutions (RE), such as self-exciting RE (SERE), which are operator-valued stochastic processes that possess two traits, self-exciting and clustering ones. These features make SERE completely different from other RE studied in the literature \cite{H,P,KS1,KS2,S00}.  The results for SERE obtained in the paper  (averaging (Theorem 1) and diffusion approximation (Theorem 2)) generalize many results in the literature devoted to the Hawkes-based models in finance and insurance, such as, e.g.,  risk model based on general compound Hawkes process \cite{S2} (self-exciting impulse traffic/transport process (8) in this paper) and general compound Hawkes processes in limit order books, see \cite{S8} (self-exciting summation on a Markov chain in this paper). We note that SERE is based on a new process $x(t)$, introduced  in the paper, which is constructed as a superposition of a Markov chain $x_n$ and a Hawkes process $N(t),$ i.e., $x(t):=x_{N(t)}.$ We call this process self-walking imbedded semi-Hawkes process (Swish Process or SwishP). 

The Random Evolutions (REs) can be described by two objects: 1) operator dynamical system $V(t)$ \index{operator dynamical system} and 2) random process $x(t).$ Depending on structure of $V(t)$ and properties of the stochastic process $x_t$ we have different kind of REs: continuous, discrete, Markov, semi-Markov, etc. 

From mathematical point of view, a RE is a solution of stochastic operator integral equation in a Banach space. The operator coefficients of such equations depend on random parameters. The RE, from physical point of view, is a model for a dynamical system whose state of evolution is subject to random variations. 
Such systems arise in many branches of science, e.g., random Hamiltonian and Shr\"{o}dinger's equations with random potential in quantum mechanics, Maxwell's equation with a random reflective index in electrodynamics, transport equation, storage equation, etc. There is a lot of applications of REs in financial and insurance mathematics \cite{S00}. One of the recent applications of RE is associated with geometric Markov renewal processes which are a regime-switching models for a stock price in financial mathematics \cite{S00}. Another recent applications of RE is a semi-Markov risk process in insurance mathematics \cite{S00}, and regime-switching financial models, including such models with jumps \cite{S00}. The REs are also examples of more general mathematical objects such as multiplicative operator functionals (MOFs) \cite{P}, which are random dynamical systems in Banach space. Different REs were introduced and studied by Hersh R. \cite{H}. Semi-Markov REs were considered in \cite{KS1,KS2}.  Inhomogeneous random evolutions, limit theorems for them  and their applications were considered in \cite{S11}. 
Some results on random operators were introduced in many papers and books, for example, random linear operators were studied in \cite{Sk2}, and integration theory in Hilbert spaces were investigated in \cite{Sk3}.

This paper is devoted to the study of a new so-called self-exciting REs and their applications. We introduce a new random process $x(t)$ such that it is based on  a superposition of a Markov chain $x_n$ and Hawkes process $N(t),$ i.e., $x(t):=x_{N(t)}.$ We call this process self-walking imbedded semi-Hawkes process (Swish Process or SwishP). Then the self-exciting REs can be constructed in similar way as, e.g., semi-Markov REs, but instead of semi-Markov process $x(t)$ we have SwishP. We give first classifications and examples of self-exciting REs (SEREs). Then we consider limit theorems for SEREs such as  averaging and diffusion approximation. Applications of SEREs are devoted to the so-called self-exciting traffic/transport process and self-exciting summation on a Markov chain, which are examples of continuous and discrete SERE. From these processes we can construct many other self-exciting processes, e.g., such as impulse traffic/transport process, self- exciting risk process, etc. We present averaged and diffusion approximation for self-exciting processes. 

We note that Hawkes process was introduced in \cite{H1,H2} and was studied in many papers and books, and was applied to many problems in finance and insurance \cite{H3}. For example, risk process based on general compound Hawkes process was studied in \cite{S2}, general compound Hawkes process applied to LOBs was intensively studied in \cite{S8}. Modelling financial contagion using mutually exciting jump processes were considered in \cite{AS}. Multivariate Hawkes processes with  application to financial data were studied in \cite{EL}. Stability of nonlinear Hawkes processes were first considered in \cite{BM}. Some limit theorems for Hawkes processes and their application to financial statistics were investigated in \cite{BDHM}. Applications of limit theorems for marked Hawkes processes to a risk model was considered in \cite{K}.

The paper is organized as follows. Section 2 gives main definitions, conditions and classifications for SERE, including definition of  a new Swish process. Section 3 presents some example of SERE, such as self-exciting impulse traffic/transport process, self-exciting geometric compound Swish process, self-exciting risk process, to name a few.. Section 4 describes martingale properties and weak convergence of SERE.  Two main limit theorems for SERE, such as averaging (Theorem 1) and diffusion approximation (Theorem 2), are presented in section 5. Applications of limit theorems for SEREs are considered in section 6. These applications are devoted to the self-exciting traffic/transport process and self-exciting summation on a Markov chain, which are continuous and discrete SERE.  Section 7, Appendix, contains proofs for Theorem 1 and Theorem 2. Discussions are considered in section 8 and section 9 concludes the paper.


\section{Main Definitions and Conditions for Swish Process and SERE}

\subsection{ Self-Walking Imbedded Semi-Hawkes Process (SwishP)}

Let $(\Omega,{\cal F}, ({\cal F}_t)_{t\geq 0},{\cal P})$ be a probability space, $t\in R_{+}:=
[0,+\infty]$, and $(X,{\cal X})$ be a measurable phase space.

{\bf Definition 1 (Hawkes Process \cite{H1,H2})}. The Hawkes process $N(t)$ is a counting process with stochastic intensity function $\lambda(t):$
$$
\lambda(t)=\lambda+\int_0^t\mu(s)dN(s),
$$
where $\lambda>0$ is a background intensity and $\mu(t)$ is a self-exciting deterministic function. We assume that $\mu(t)=\alpha e^{-\beta t}, \alpha, \beta>0,$ is an exponential self-exciting function. We suppose that $\int_0^{+\infty}\mu(s)ds:=\hat\mu<1$- stability condition to avoid accumulating of jumps. Here, $\hat\mu=\alpha/\beta.$

Moments of jumps (or arrival times of events) of the Hawkes process are $\tau_n,$ and inter-arrival times are $\theta_n\in R_+,$ therefore $\tau_{n}:=\sum_{k=0}^{n}\theta_{k}.$ Let we also have homogeneous Markov chain $x_n\in X,$ where $X$ is a space of states (or phase space) for $x_n.$

{\bf Definition 2 (SwishP)}. Swish Process (SwishP) is defined in the following way:
\begin{equation}\label{e1}
x(t):=x_{N(t)},
\end{equation}
where $x_n$ is a Markov chain and $N(t)$ is a Hawkes process, and they are independent.

As long as we suppose  that $x_{n}$ and $\theta_{n}, n\geq 0,$ are independent, then

\begin{equation}\label{e2}
Q(x,A,t):={\cal P}\{x_n\in A, \theta_n\leq t|x_0=x\}=P(x,A)G(t),
\end{equation}
where
\begin{equation}\label{e3}
\begin{array}{rcl}
P(x,A)&:=&{\cal P}(x_{n+1}\in A/x_{n}=x), A\in {\cal X},\\
G(t)&:=&{\cal P}(\theta_{n+1}\leq t).
\end{array}
\end{equation}

\subsection{Definition and Classification of SERE}

Let  $(B,\cal B, \|\cdot\|)$ be a separable Banach space \cite{A,DS, S11}.

Let $\{\Gamma(x);x\in X\}$ be a family of operators \index{family of operators} on the dense subspace 
$B_{0}\in B$, which is common domain for $\Gamma(x)$, independent of $x$,
non-commuting and unbounded in general, such that map $\Gamma(x)f:X\to B$ is
strongly ${\cal X}/\cal B$-measurable for all $f\in B_0$, $\forall t\in R_{+}.$ Also,
let $\{{\cal D}(x);x\in X\}$ is a family of bounded linear operators on $B$, 
such that map ${\cal D}(x)f:X\to B$ is ${\cal X}/\cal B$-measurable,
$\forall f\in B.$

{\bf Definition 3 (Self-Exciting Random Evolution (SERE)}. SERE are defined by the solution of the following stochastic operator
integral equation in separable Banach space $B$:
\begin{equation}\label{e4}
V(t)f=f+\int_{0}^{t}V(s)\Gamma(x(s))fds+
\sum_{k=1}^{N(t)}V(\tau_{k})[{\cal D}(x_{k})-I]f,
\end{equation}
where $N(t)$ is a Hawkes process, $x(t)$ is a Swish process, $x_k$ is a Markov chain, $I$ is an identity operator on $B$, $f\in B_0.$ Because of the stability condition $\int_0^{+\infty}\mu(s)ds:=\hat\mu<1,$ the sum in (\ref{e4}) is finite.


If $x(t)$ in (\ref{e1}) is a Markov or semi-Markov process, then SERE in 
(\ref{e4}) is called a Markov or semi-Markov RE, respectively. In this case, $N(t)$ is a renewal process, \cite{S00}.
If ${\cal D}(x)\equiv I$, $\forall x\in X$, then $V(t)$ in (\ref{e4}) is
called a continuous SERE. If $\Gamma(x)\equiv 0$, $\forall x\in X,$ is a zero operator on $B$, then $V(t)$
in (\ref{e4}) is called a jump RE. SERE $V_{n}:=V(\tau_{n})$ is called a discrete SERE.
Operators $\Gamma(x)$, $x\in X$, describe a continuous component $V^{c}(t)$ of
SERE $(V(t))_{t\geq 0}$ in (\ref{e4}), and operators ${\cal D}(x)$ describe a jump component
$V^{d}(t)$ of SERE $(V^{d}(t))_{t\geq 0}$ in (\ref{e4}).

In such a way, SERE is described by two objects:
1) operator dynamical system $(V(t))_{t\geq 0}$ \index{operator dynamical system};
2) random process $x(t)$.

We note, that it turned out to be (see \cite{KS1} for more details) 
\begin{equation}\label{e5}
V(t)=\Gamma_{x(t)}(t-\tau_{N(t)})
\prod_{k=1}^{N(t)}{\cal D}(x_{k})\Gamma_{x_{k-1}}(\theta_{k}),
\end{equation}
where $\Gamma_{x}(t)$ are the semigroups of operators of $t$ generated by the 
operators $\Gamma(x)$, $\forall x\in X.$ Because of the stability condition $\int_0^{+\infty}\mu(s)ds:=\hat\mu<1,$ the product in (\ref{e5}) is finite.

We also note, that SERE in (\ref{e5}) is called a discontinuous SERE.
Under conditions introduced in the next section, the solution $V(t)$ of the equation (\ref{e4})
is unique and can be represented by product (\ref{e5}), that can be proved by
constructive or iterative methods, similarly to the approach which was done in \cite{KS1, KS2} for semi-Markov REs.


\subsection{Conditions for SERE}

We suppose that the following conditions are satisfied (which will be used in Theorem 1 (Averaging of SERE) and Theorem 2 (DA of SERE)):

A) there exists Hilbert spaces $H$ and $H^{*}$ such that compactly
imbedded in Banach spaces $B$ and $B^{*}$ respectively,
$H \subset B, \quad H^{*} \subset B^{*}$, where $B^{*}$ is a dual space
to $B$, that divides points of $B$;

B) operators $\Gamma(x)$ and$(\Gamma(x))^{*}$ are dissipative on any
Hilbert space $H$ and $H^{*}$ respectively;

C) operators ${\cal D}(x)$ and ${\cal D}^{*}(x)$ are contractive on any
Hilbert space $H$ and $H^{*}$ respectively;

D) $(x_{n}; \quad n \geq 0)$ is a uniformly ergodic Markov chain
with stationary distribution $\rho (A),\quad A \in \cal{X}$;

E) $m_{i}: = \int_{0}^{\infty} t^{i} G(dt)$ are uniformly
integrable, $\forall i = 1, 2, 3$, where

\begin{eqnarray}\label{e29}
G(t): = {\cal P}\{\omega: \theta_{n+1} \leq t\};
\end{eqnarray}

F) \begin{eqnarray}\label{e30}
\int_{X}\rho(dx)\|\Gamma(x)f\|^{k} < +\infty;
\int_{X}\rho(dx)\|P{\cal D}_{j}(x)f\|^{k} < +\infty;\quad j=1,2,\nonumber\\
\int_{X}\rho(dx)\|\Gamma(x)f\|^{k-1}\cdot\|P{\cal D}_{j}(x)f\|^{k-1}
< +\infty;\quad \forall k = 1, 2, 3, 4, f \in B_0,
\end{eqnarray}
where $P$ is on operator generated by the transition probabilities
$P(x, A)$ of Markov chain $(x_{n}; \quad n \geq 0)$:

\begin{eqnarray}\label{e31}
P(x, A): = {\cal P}\{\omega: x_{n+1} \in A/x_{n} = x\},
\end{eqnarray}
and ${\cal D}_{j}(x); \quad x \in X, \quad j = 1, 2\}$ is a family
of some closed operators;

G) $\int_0^{+\infty}\mu(s)ds<1,$ where $\mu(t)$ is a self-exciting exponential function for Hawkes process $N(t);$

H) $\int_0^{+\infty}s\mu(s)ds<+\infty.$

If $B: = C_{0}(R)$, then $H: = W^{l, 2}(R)$ is a Sobolev space \cite{A,DS,So}, and
$W^{l, 2}(R) \subset C_{0}(R)$ and this imbedding is compact.
For the spaces $B: = L_{2}(R)$ and $H: = W^{l, 2}(R)$ we have the
same situation.

Let us explain why we need all the conditions A)-H).
It follows from the conditions A) - B) that operators $\Gamma(x)$ and
$(\Gamma(x))^{*}$ generate a strongly continuous contractive
semigroup of operators $\Gamma_{x}(t)$ and $\Gamma_{x}^{*}(t), \quad
\forall x \in X$, in $H$ and $H^{*}$ respectively. From the
conditions A)--C) it follows that SERE $V(t)$ in (1) is a contractive
operator in $H, \quad \forall t \in R_{+}$, and $\|V(t)f\|_{H}$ is
a semimartingale $\forall f \in H$. In such a way, the conditions
A) - C) supply the following result:

{\bf Proposition (Tightness of SERE).} SERE $V(t) f$ is a tight process\index{tight process} in $B$, namely, $\forall \Delta>0$
there exists a compact set\index{compact set} $K_{\Delta}$:

\begin{eqnarray} \label{e32}
{\cal P}\{V(t)f \in K_{\Delta}; 0 \leq t \leq T\}\geq 1 - \Delta.
\end{eqnarray}

{\bf Proof.} This result follows from Kolmogorov - Doob inequality for 
semimartingale $\|V(t) f\|_{H}$ \cite{KS1}.

Condition (\ref{e32}) is one of the main steps in the proving of limit theorems
and rates of convergence for the sequence of SERE in series scheme.
Conditions D)-E) need for asymptotical results for SERE associated with Markov chain $x_n$ and r.v. $\theta_n.$ We require technical condition D) for the averaged by $\rho(dx)$ moments  for the operators $\Gamma(x)$ and ${\cal D}_j(x), j=,1,2,$ in the asymptotic expansions of operators $\Gamma_x(t)$ (see (26)) and ${\cal D}^{\epsilon}(x)$ (see (28) and (32)). Technical conditions G) and H), designed for the Hawkes process $N(t),$ needed for the asymptotic analysis of HP in the averaging (Theorem 1) and DA (Theorem 2) results for SERE.

\section{Some Examples of SERE}

Many examples for SERE can be obtained by assigning to every operator $\Gamma(x)$ and ${\cal D}(x)$ a specific map. For example, if we take
$$\Gamma(x):=v(x)\frac{d}{dz}, \quad {\cal D}(x)\equiv I, \quad B=C^{1}(R),$$
where $C^1(R)$ is a space of differentiable functions on $R,$ then we get a transport equation for (4) describing a motion of a particle with random velocity $v(x(t)),$ i.e.,
$$z(t)=z+\int_{0}^{t}v(z(s),x(s))ds.$$ Equation (4) in this case is
$$
f(z(t))=f(z)+\int_{0}^{t}v(z(s),x(s))\frac{d}{dz}f(z(s))ds.
$$
Below we present many other examples of SERE.

{\bf Example 1. Self-exciting Impulse Traffic/Transport process.}
Let $B=C(R)$ be a space of continuous functions on $R,$ and operators $\Gamma(x)$ and ${\cal D}(x)$ are defined by the 
following way:
\begin{equation}\label{e6}
\Gamma(x)f(z):=v(z,x)\frac{d}{dz}f(z),\quad 
{\cal D}(x)f(z):=f(z+a(x)),
\end{equation}
where functions $v(z,x)$ and $a(x)$ are continuous and bounded on $R\times X$
and $X$ respectively, $\forall z\in R$, $\forall x\in X$, $ f(z)\in
C^{1}(R):=B_{0}.$
Then the equation (\ref{e4}) takes the form:

\begin{equation}\label{e7}
f(z(t))=f(z)+\int_{0}^{t}v(z(s),x(s))\frac{d}{dz}f(z(s))ds+
\sum_{k=1}^{N(t)}[f(z(\tau_{k})+a(x_{k}))-f(z(\tau_{k}))],
\end{equation}
and SERE $V(t)$ is defined by the relation:
$$V(t)f(z)=f(z(t)),\quad z(0)=z.$$
Equation (\ref{e7}) is a functional form for self-exciting impulse traffic process 
$z_{t},$ which satisfies the equation:

\begin{equation}\label{e8}
z(t)=z+\int_{0}^{t}v(z(s),x(s))ds+\sum_{k=1}^{N(t)}a(x_{k}).
\end{equation}
We note that self-exciting impulse traffic/transport process $z(t)$ in (12) is a realization of
discontinuous RE. We call $z(t)$ a self-exciting traffic/transport process whenever $a(x)=0.$ We note, that if we take $v(z,x)=v$ and $a(x)$ takes negative values, then we get risk model based on general compound Hawkes process \cite{S2}. We note that sum in (12) cannot be exploded because of the stationary condition $\hat\mu<1$ (see condition G)).

{\bf Example 2. Self-Exciting Summation on a Markov chain.} \index{Summation on a Markov chain}
Put $v(z,x)\equiv 0$, $\forall z\in R$, $\forall x\in X$, in (12).
Then the process
\begin{equation}\label{e9}
z(t)=z+\sum_{k=1}^{N(t)}a(x_{k})
\end{equation}
is a summation on a Markov chain $(x_{n};n\geq 0)$ and it is a realization of a
jump SERE.
We also called it general compound Hawkes process \cite{S8}.

{\bf Example 3. Self-Exciting Diffusion Process.}
Let $B=C(R)$, $B_{0}=C^{2}(R)$-twice differentiable functions on $R,$ $P_{x}(t,z,A)$ be a distribution
function or transition probability kernel, which respects to the diffusion process $\xi (t)$, that is the solution
of the stochastic differential equation in $R$ with SwishP switchings:
\begin{equation}\label{e10}
d\xi (t)=a(\xi(t),x(t))dt+b(\xi(t),x(t))dw_{t}, \xi(0)=z,
\end{equation}
where $x(t)$ is a Swish process independent of a standard Wiener process
$w_{t}$, coefficients $a(z,x)$ and $b(z,x)$ are bounded and continuous
functions on $R\times X$.
Let us define the following contraction semigroups of operators on $B$:
\begin{equation}\label{e11}
\Gamma_{x}(t)f(z):=\int_{R}P_{x}(t,z,dy)f(y), f(y)\in B, x\in X.
\end{equation}
Their infinitesimal operators $\Gamma(x)$ have the following form:
$$\Gamma(x)f(z)=a(z,x)\frac{d}{dz}f(z)+
2^{-1}b^{2}(z,x)\frac{d^{2}}{dz^{2}}f(z), f(z)\in B_{0}.$$
Here, the operators 
${\cal D}(x)\equiv I$, $\forall x\in X$, are identity operators.
Then the equation (\ref{e4}) takes the form:
\begin{equation}\label{e12}
f(\xi(t))=f(z)+\int_{0}^{t}[a(\xi(s),x(s))\frac{d}{dz}+
2^{-1}b^{2}(\xi(s),x(s))\frac{d^{2}}{dz^{2}}]f(\xi(s))ds,
\end{equation}
and SERE $V(t)$ is defined by the relation
$$V(t)f(z)=E[f(\xi(t))/x(s);0\leq s\leq t; \xi(0)=z].$$
Equation (\ref{e12}) is a functional form for diffusion process $\xi(t)$ in
(\ref{e10}) in SwishP random media $x(t).$
We note that diffusion process $\xi(t)$ in (\ref{e10}) is a realization of
continuous SERE.

{\bf Example 4. Self-Exciting Geometric Compound Process (SEGCP)}.
Let $x(t)$ be a Swish process on the phase space $X\times \bf{R}^+,$  i.e., $x(t):=x_{Nt)}.$ Let $c(x)$ be a bounded continuous function on $X$ such that $c(x)>-1.$ We define a SEGCP $S_t$ with SwishP $x(t)$ as follows:
\begin{equation}\label{e13}
S_t:=S_0\prod_{k=0}^{N(t)}(1+c(x_k)),
\end{equation}
where $S_0>0$ is the initial value of $S_t.$ We call  the process Self-Exciting Geometric Compound Process (SEGCP) which is similar to the geometric compound Poisson process
$$
S_t=S_0\prod_{k=1}^{N(t)}(1+Y_k),
$$
where $S_0>0,$ $N(t)$ is a standard Poisson process, $(Y_k)_{k\in Z_{+}}$ are i. i. d. random variable,
which is a trading model in many financial applications as a pure jump model (see \cite{S00}).  We note that process $S_t$ cannot be exploded because of the stationary condition $\hat\mu<1.$

Let $B:=C_0(\bf{R}_+)$ be a space of continuous functions on $\bf{R}_+,$ vanishing at  infinity, and let us define
a family of bounded contracting operators $D(x)$ on $C_0(\bf{R}_+):$

\begin{equation} \label{e14}
D(x)f(s):=f(s(1+c(x))),~ x\in X, s\in \bf{R}_+.
\end{equation}

With these contraction operators $D(x)$ we define the following jump SERE $V(t)$ of SEGCP in (\ref{e13})
\begin{equation}\label{e15}
V(t)=\prod_{k=0}^{N(t)}D(x_k):=D(x_{N(t)})\circ D(x_{N(t)-1})\circ\ldots \circ D(x_1)\circ D(x_0).
\end{equation}
Using  (\ref{e14}) we obtain from (\ref{e15})
\begin{equation} \label{e16}
V(t)f(s)=\prod_{k=0}^{N(t)}D(x_k)f(s)=f(s\prod_{k=0}^{N(t)}(1+c(x_k))=f(S_t),
\end{equation}
where $S_t$ is defined in (\ref{e13}) and $S_0=s.$

{\bf Example 5. Self-Exciting Risk Process (SERP)}. SERSP is defined as 

\begin{equation}\label{e8}
R(t)=u+ct-\sum_{k=1}^{N(t)}a(x_{k}),
\end{equation}
where $u$ is an initial capital of an insurance company, $c$ is a premium, $a(x)$ is a measurable bounded function on $X,$ and $N(t)$ is a hakes process. Here,
$$
\Gamma_x(t)f(u)=\Gamma(t) f(u):=f(u+ct), \quad {\cal D}(x)f(u):=f(u-a(x)).
$$


\section{Martingale Property and Weak Convergence of SERE}


We use martingale methods to study SERE and limit theorems for them.
The main idea is that process 
\begin{equation}\label{e17}
M_{n}f:=V_{n}f-f-\sum_{k=0}^{n-1}E[V_{k+1}f-V_{k}f/{\cal F}_{k}], V_{0}=I, f\in B,
\end{equation}
is an ${\cal F}_{n}$-martingale in $B$, where $\sigma$-algebra
$${\cal F}_{n}:=\sigma\{x_{k},\tau_{k};0\leq k\leq n\}, V_{n}:=V(\tau_{n}),$$
$E$ is an expectation by probability $\cal P.$  It is defined as a Bochner integral in Banach space $B$ (see \cite{S11}). For definitions of expectations of operator-valued processes and martingale meaning in this context, see \cite{S11}, Chapter 1, sec.1.1
Representation of the martingale $M_{n}$ (see (22)) in the form of
martingale-difference
\begin{equation}\label{e18}
M_{n}f=\sum_{k=0}^{n-1}[V_{k+1}f-E(V_{k+1}f/{\cal F}_{k})]
\end{equation}
gives us the possibility to calculate the weak quadratic variation:
\begin{equation}\label{e19}
<l(M_{n}f)>:=\sum_{k=0}^{n-1}E[l^{2}((V_{k+1}-V_{k})f)/{\cal F}_{k}],
\end{equation}
where $l\in B^{*}$ is a linear functional from $B^*,$ and $B^{*}$ is a dual space to $B$, dividing points of $B.$

In general, the weak convergence of SERE in series scheme is obtained from the criterion of weak
compactness of the processes with values in separable Banach space.
The limiting SERE we obtain from the solution of a martingale problem in the form of
some integral operator equations in Banach space $B$.
In particular, the martingale method of obtaining of the limit theorems for the sequence of SERE
is based on the solution of the following problems (see \cite{S11}, Chapter 1, sec. 1.2):
1) weak compactness of the family of measures generated by the sequences of SERE (see \cite{S11}, Chapter 1, Theorem 1.9, and {\bf Proposition (Tightness of SERE)}, sec. 2.3);
2) any limiting point of this family of measures is the solution of a martingale
problem;
3) the solution of martingale problem is unique.
The conditions 1)-2) guarantee the existence of weakly convergent subsequence,
and condition 3) gives the uniqueness of the weak limit.
It follows from 1)-3) that consequence of SERE converges weakly to the unique
solution of a martingale problem.

We also use the representation
\begin{equation}\label{e20}
V_{k+1}f-V_{k}f=[\Gamma_{x_{k}}(\theta_{k+1}){\cal D}(x_{k})-I]V_{k}f,\quad
V_{k}:=V(\tau_{k}),
\end{equation}
and the following expression for semigroups of operators $\Gamma_{x}(t)$ (see \cite{DS}, \cite{S11}):
\begin{equation}\label{e21}
\Gamma_{x}(t)f=f+\sum_{k=1}^{n-1}\frac{t^{k}}{k!}\Gamma^{k}(x)f+
(n-1)^{-1}\int_{0}^{t}(t-s)^{n}\Gamma_{x}(s)\Gamma^{n}(x)fds,
\forall x\in X,
\end{equation}
 $\forall f\in \cap_{x\in X}Dom(\Gamma^{n}(x)).$
Taking into account (\ref{e17})-(\ref{e21}) and  conditions A)-H) from sec. 2, we obtain the limit theorems for SERE in the next sections.


\section{Main Results: Averaging and Diffusion Approximation of SERE}

\subsection{Averaging of Self-Exciting Random Evolutions}

Let's consider a SERE in series scheme:

\begin{eqnarray} \label{e33}
V_{\varepsilon}(t) = f + \int_{0}^{t}\Gamma
(x(s/\varepsilon)) V_{\varepsilon}(s)
fds + \sum_{k=1}^{N(t/\varepsilon)}
[{\cal D}^{\varepsilon}( x_{k}) - I] V_{\varepsilon}(\varepsilon \tau_{k})f,
\end{eqnarray}

where

\begin{eqnarray}\label{e34}
{\cal D}^{\varepsilon}(x) = I + \varepsilon D_{1}(x) + o(\varepsilon),
\end{eqnarray}
$\{D_{1}(x); x, y \in X\}$ is a family of closed linear operators,
$\|0  (\varepsilon) f\|/\varepsilon \to 0 \quad \varepsilon \to 0, \quad 
\varepsilon$ is a small parameter,

\begin{eqnarray}\label{e35}
f \in B_{0}: = \bigcap_{x, y \in X} Dom
(\Gamma^{2}(x)) \cap Dom (D_{1}^{2}(x, y)).
\end{eqnarray}

Another form for $V_{\varepsilon}(t)$ in (\ref{e35}) is (see also (5)):

\begin{eqnarray}\label{e36}
V_{\varepsilon}(t) = \Gamma_{x(t/\varepsilon)}
(t - \varepsilon \tau_{N(t/\varepsilon)})
\prod_{k=1}^{N(t/\varepsilon)} {\cal D}^{\varepsilon}
(x_{k}) \Gamma_{k-1}(\varepsilon \theta_{k}).
\end{eqnarray}

{\bf Theorem 1 (Averaging of SERE)}. Let the conditions A)-D), E) with $i=2,$ F) with $k=2, j=1,$ and G) are satisfied (see sec. 2.3). Then SERE $V_{\epsilon}(t)$ in (27) or in  (30) converges weakly as $\epsilon\to 0$ to the following limiting process $V_0(t)$ satisfying the equation

\begin{eqnarray}
 V_{0}(t)f - f - \hat\lambda\int_{0}^{t}
(\hat{\Gamma} + \hat{D})V_{0}(s)fds=0,
\end{eqnarray}
where
$$ 
\begin{array}{rcl}
\hat{\Gamma}&: =&m \int_{x}\rho(dx)\Gamma(x),\\
\hat{D}&: =& \int_{x}\rho(dx)D_{1}(x),\nonumber,\\
\end{array}
$$
$m: = m_{1},$ $\hat\lambda$ is defined by $$\hat\lambda:=\frac{\lambda}{1-\hat\mu},$$ where $\lambda$ and $\hat\mu$ are defined in sec. 2.1. Proof of Theorem 1 is given in the Appendix.

\subsection{Diffusion Approximation (Functional Central Limit Theorem) of Self-Exciting Random Evolutions}

Let us consider SERE $V^{\varepsilon}(t)$, which is defined as 
as
$$
V^{\varepsilon}(t) = f + \int_{0}^{t}\Gamma
(x(s/\varepsilon^2)) V^{\varepsilon}(s)
fds + \sum_{k=1}^{N(t/\varepsilon^2)}
[{\cal D}^{\varepsilon}( x_{k}) - I] V^{\varepsilon}(\varepsilon \tau_{k})f,
$$


\begin{eqnarray}\label{e43}
{\cal D}^{\varepsilon}(x): = I + \varepsilon D_{1}(x)+
\varepsilon^{2}D_{2}(x) + o(\varepsilon^{2}),
\end{eqnarray}
$\{D_i(x); \quad x \in X, \quad i = 1, 2\}$ are closed
linear operators and\\
$\|0(\varepsilon^{2}) f\|/\varepsilon^{2} \to 0,  \varepsilon \to 0$

\begin{eqnarray}\label{e44}
\forall f \in B_{0}: = \bigcap_{x \in X} Dom
(\Gamma^{4}(x))\bigcap Dom(D_{2}(x)),\nonumber\\
Dom(D_{2}(x)) \subseteq Dom (D_{1}(x));\quad
D_{1}(x) \subseteq Dom(D_{1}(x)),\nonumber\\
\forall x \in X, \quad \Gamma^{i}(x) \subset
Dom (D_{2}(x)), \quad i = \overline{1, 3}.
\end{eqnarray}

In such a way, $V^{\varepsilon}(t)$ can be also presented inthe following way (see also (5)):

\begin{eqnarray}\label{e45}
V^{\varepsilon}(t) = \Gamma_{x(t/\varepsilon^{2})}
(t/\varepsilon - \varepsilon \tau_{N(t/\varepsilon^{2})})
\prod_{k=1}^{N(t/\varepsilon^{2})} {\cal D}^{\varepsilon}
(x_{k}) \Gamma_{x_{k-1}}(\varepsilon\theta_{k}),
\end{eqnarray}
where ${\cal D}^{\varepsilon}(x)$ are defined in (32).

{\bf Theorem 2 (Diffusion Approximation of SERE)}. Let the conditions A)-D), E) with $i=3,$ F) with $k=3, j=2,$ and G)-H) are satisfied (see sec. 2.3). Also, let the {\bf balance condition} be satisfied:

\begin{eqnarray}\label{e46}
\int_{X}\rho(dx)[m\Gamma(x) + D_{1}(x)] f =0,\quad
\forall f \in B_{0}
\end{eqnarray}

Then the SERE $V^{\varepsilon}(t)$ converges weakly as $\epsilon\to 0$ to the limiting process $V^0(t)$ which satisfies the following equation 

\begin{eqnarray} 
V^{0}(t)f = f + \hat\lambda\cdot \int_{0}^{t}\hat{L}\cdot V^{0}(s)
fds + \int_{0}^{t}\int_{X}\sigma\Gamma(x)V^{0}(s)f W(dx, ds),
\end{eqnarray}

where 
$$
\begin{array}{rcl}
\hat L: &=& \int_{X}\rho(dx)L(x),\\
L(x) &:=& (m\Gamma(x) + D_{1}(x))
(R_{0} - I)(m\Gamma(x) + D_{1}(x)) +\nonumber\\ 
&\mbox{}& + m_{2}
\Gamma^{2}(x)/2 + mD_{1}(x)
\Gamma(x) + D_{2}(x),\nonumber
\end{array}
$$
$R_{0}$ is a potential operator of $(x_{n};\quad n \geq 0)$ (see \cite{KS1,KS2,S00}),
$W(dx, ds)$ is the Wiener orthogonal martingale measure with quadratic variation $\rho(dx)ds$ (see \cite{KS2,S11}) and $\hat\lambda$ is defined in (31). Proof of Theorem 2 is given in the Appendix.

\section{Applications of SERE in Limit Theorems}

Applications of SERE will be given to self-exciting traffic/transport process and self-exciting summation on a Markov chain, which are continuous and discrete SERE.
From these processes we can construct many other self-exciting processes, e.g., such as impulse traffic/transport process, self-exciting risk process, etc.
We will present averaged and diffusion approximation of self-exciting processes. All convergences are in Skorokhod topology \cite{Sk} when $\epsilon\to 0.$
The results here follow from sec. 5.1-5.2 (Theorem 1 and Theorem 2). 
\subsection{Averaged Self-exciting Traffic/Transport and Summation on a Markov Chain Processes (see sec. 4.2)}

If we have self-exciting traffic/transport process in series scheme $t/\epsilon$
\begin{equation}\label{e79}
z^{\epsilon}(t)=z+\int_{0}^{t}v(z^{\epsilon}(s),x(s/\epsilon))ds,
\end{equation}
then the averaged self-exciting traffic/transport process satisfies the ODE (weak convergence when $\epsilon\to 0$ in Skorokhod topology \cite{Sk})

\begin{equation}\label{e80}
\frac{d\hat z(t)}{dt}=\hat v(\hat z(t)),\quad \hat z(0)=z,
\end{equation}
where $\hat v(z):=\hat\lambda m\int_X\rho(dx)v(z,x),$ $\hat\lambda:=\frac{\lambda}{1-\hat\mu}.$

If we have self-exciting summation on a Markov chain process in series scheme $t/\epsilon,$

\begin{equation}\label{e81}
z^{\epsilon}_{t}=z+\epsilon\sum_{k=1}^{N(t/\epsilon)}a(x_{k}),
\end{equation}
then the averaged self-exciting summation on a Markov chain process has the following form: (weak convergence when $\epsilon\to 0$ in Skorokhod topology \cite{Sk})

\begin{equation}\label{e82}
\hat z_t=z+\hat at,
\end{equation}
where $\hat a:=\hat\lambda\int_X\rho(dx)a(x).$

\subsection{Diffusion Self-exciting Traffic/Transport and Summation on a Markov Chain Processes}

If we have self-exciting traffic/transport process in series scheme $t/\epsilon^2$

\begin{equation}\label{e83}
z^{\epsilon}(t)=z+\int_{0}^{t}v(z^{\epsilon}(s),x(s/\epsilon^2))ds,
\end{equation}

with balance condition $\hat v(z)=\int_X\rho(dx)v(z,x)=0,$ then the diffusion self-exciting traffic/transport process has the following SDE form (weak convergence when $\epsilon\to 0$ in Skorokhod topology \cite{Sk})

\begin{equation}\label{e84}
d\hat z(t)=b(\hat z(t))dt+\sigma(\hat z(t))dw(t),
\end{equation}

where $b(z):=\hat\lambda \int_X\rho(dx)[m^2v(z,x)(R_0-I)v'(z,x)+m_2v(z,x)v'(z,x)/2]$ and $\sigma^2(z):=2\hat\lambda \int_X\rho(dx)[m^2v(z,x)(R_0-I)v(z,x)+m_2v^2(z,x)/2],$ $R_0$ is a potential of the Markov chain $(x_n, n\geq 0),$ $w(t)$ is a standard Wiener process.  

If we have  self-exciting summation on a Markov Chain Processes in series scheme $t/\epsilon^2$

\begin{equation}\label{e85}
z^{\epsilon}_{t}=z+\epsilon\sum_{k=1}^{N(t/\epsilon^2)}a(x_{k}),
\end{equation}
with a balance condition $\hat a=\int_X\rho(dx)a(x)=0,$ then the diffusion summation on a Markov chain  (weak convergence when $\epsilon\to 0$ in Skorokhod topology \cite{Sk}) is a Wiener process with variance

\begin{equation}\label{e86}
\sigma^2:=2\hat\lambda\int_X\rho(dx)[m^2a(x)(R_0-I)a(x)+m_2a^2(x)/2],
\end{equation}
i.e., $z^{\epsilon}_{t}\to \hat z_t$ when $\epsilon\to 0,$ and $$\hat z_t=z+\sigma w(t). $$

\section{Appendix: Proofs of Theorem 1 (Averaging of SERE) and Theorem 2 (DA of SERE)}


\subsection{Proof of Averaging of SERE (Theorem 1)}

Under conditions A) - C) the sequence of SERE $V_{\varepsilon}(t)f$ is
tight (see {\bf Proposition (Tightness of SERE)}, sec. 4). 

Under conditions $D), E), i = 2, F), k=2, j = 1$, the sequence of
SERE $V_{\varepsilon}(t)f$ is weakly compact $\rho - a. s.$ in
$D_{B}[0, + \infty)$ with limit points in\\
$C_{B}[0, +\infty),\quad f \in B_{0}$. Here: $D_{B}[0, + \infty)$ is a space of $B$-valued c\`{a}dl\`{a}g function on $[0.+\infty),$ and  $C_{B}[0, +\infty)$ is a space of $B$-valued continuous function on $[0,+\infty).$

Let's consider the following process in $D_{B}[0, +\infty)$: 

\begin{eqnarray}\label{e37}
M_{\nu(t/\varepsilon)}^{\varepsilon} f^{\varepsilon}: =
V_{\nu(t/\varepsilon)}^{\varepsilon} f^{\varepsilon} -
f^{\varepsilon} - \sum_{k=0}^{N(t/\varepsilon) - 1}
E_{\rho}[V_{k+1}^{\varepsilon} f_{k+1}^{\varepsilon} - 
V_{k}^{\varepsilon} f_{k}^{\varepsilon}/{\cal F}_{k}],
\end{eqnarray}
where $V_{n}^{\varepsilon}: = V_{\varepsilon}(\varepsilon \tau_{n}),$ 
$$f^{\varepsilon}: = f + \varepsilon f_{1}(x(t/\varepsilon)),$$
$$f_{k}^{\varepsilon}: = f^{\varepsilon}(x_{k}),$$ function $f_{1}(x)$ is
defined from the equation

\begin{eqnarray}\label{e38}
(P - I)f_{1}(x) = 
\left[
(\hat{\Gamma} + \hat{D}) - (m\Gamma(x) + PD_{1}(x))
\right]f,\nonumber\\
\hat{\Gamma}: =m \int_{x}\rho(dx)\Gamma(x),\quad 
\hat{D}: = \int_{x}\rho(dx)D_{1}(x),\nonumber\\
m: = m_{1}
\end{eqnarray}
(see E)), $f \in B_{0}$.

The process $M_{N(t/\varepsilon)}^{\varepsilon} f^{\varepsilon}$ is
an ${\cal F}_{t}^{\varepsilon}$--martingale with respect to the
$\sigma$--algebra ${\cal F}_{t}^{\varepsilon} := \sigma
\{x(s/\varepsilon); 0 \leq s \leq t\}$.

The martingale $M_{N(t/\varepsilon)}^{\varepsilon} f^{\varepsilon}$
in (\ref{e37}) has the asymptotic representation:

\begin{eqnarray}\label{e39}
M_{N(t/\varepsilon)}^{\varepsilon} f^{\varepsilon} =
V_{N(t/\varepsilon)}^{\varepsilon} f - f - \varepsilon
\sum_{k=0}^{N(t/\varepsilon)} (\hat{\Gamma} +\hat{D})
V_{k}^{\varepsilon} f + 0_{f}(\varepsilon),
\end{eqnarray}
where $\hat{\Gamma}, \hat{D}, f, f^{\varepsilon}$ are defined in
(\ref{e38}) and $$\|0_{f}(\varepsilon)\| / \varepsilon \to const\quad
as \varepsilon \to 0, \quad \forall f \in B_{0}.$$


The families $l(M_{N(t/\varepsilon)}^{\varepsilon}
f^{\varepsilon})$ and 
$$
l
\left(
\sum_{k=0}^{N(t/\varepsilon)} E_{\rho}
[(V_{k+1}^{\varepsilon} f_{k+1}^{\varepsilon} - 
V_{k}^{\varepsilon} f_{k}^{\varepsilon})/{\cal F}_{k}]
\right)
$$
are weakly compact for all $l \in B_{0}^{*}$ is a some dense subset
from $B^{*}$.
Let $V_{0}(t)$ be a limit process for $V_{\varepsilon}(t) as\quad\varepsilon \to 0$.

Since (see (25))

\begin{eqnarray}\label{e40}
[V_{\varepsilon} (t) - V_{N(t/\varepsilon)}^{\varepsilon}] = 
[\Gamma_{x(t/\varepsilon)} (t - \varepsilon\tau_{N(t/\varepsilon)})
- I] \cdot V_{N(t/\varepsilon)}^{\varepsilon}
\end{eqnarray}
and the right-hand side in (\ref{e40}) tends to zero $as\quad \varepsilon \to
0$, then it is clearly that the limits for $V_{\varepsilon} (t)$ and
$V_{N(t/\varepsilon)}^{\varepsilon}$ are the same, namely,
$V_{0}(t)\quad \rho - a. s.$

Taking into account the LLN for Hawkes process, $N(t)/t\to_{t\to+\infty}\lambda/(1-\hat\mu),$ see, e.g., \cite{BDHM, BM,FLM,Z}, the sum $\varepsilon \cdot \sum_{k=0}^{N(t/\varepsilon)}
(\hat{\Gamma} + \hat{D}) V_{k}^{\varepsilon} f$ converges
strongly $as\quad \varepsilon \to 0$ to the integral
$$
\hat\lambda\cdot \int_{0}^{t}(\hat{\Gamma} + \hat{D})
V_{0}(s)fds,
$$ 
where
\begin{eqnarray}\label{e41}
\hat\lambda:=\frac{\lambda}{1-\hat\mu}.
\end{eqnarray}

The quadratic variation of the martingale 
$l(M_{\nu(t/\varepsilon)}^{\varepsilon} f^{\varepsilon})$ tends to
zero, and, hence, 

$$
M_{\nu(t/\varepsilon)}^{\varepsilon} f^{\varepsilon} \to 0\quad
as\quad \varepsilon \to 0, \quad \forall f \in B_{0}, \quad
\forall e \in B_{0}^{*}.
$$

Passing to the limit in (\ref{e39}) as $\varepsilon \to 0$ and taking
into account the all previous reasonings we obtain that the limit
process $V_{0}(t)$ satisfies the equation:

\begin{eqnarray}\label{e42}
V_{0}(t)f - f - \hat\lambda\int_{0}^{t}
(\hat{\Gamma} + \hat{D})V_{0}(s)fds=0,
\end{eqnarray}
where $\hat\lambda$ is defined in (\ref{e41}).

The uniqueness of the limiting SERE $V_{0}(t)f$ in 
averaging
scheme follows from the equation (\ref{e42}) and the fact that if the 
operator $\hat{\Gamma} + \hat{D}$ generates a semigroup,
then $V_{0}(t)f = \exp\{\hat\lambda(m\hat{\Gamma} + \hat{D}))\cdot t\} f$ and
this representation is unique. {\it Q.E.D.}

{\bf Remark.} In a similar way we can consider averaging of self-exciting random evolutions in reducible phase space, or
merged self-exciting random evolutions, as it was done in \cite{KS1, KS2,S00} for semi-Markov RE.

\subsection{Proof of DA of SERE (Theorem 2)}

 Under conditions A) - C) the sequence of SERE $V^{\varepsilon}
(tf$ is tight (see {\bf Proposition (Tightness of SERE)}, sec. 4).

Under conditions $D), E), i = 3, F), k=4, j=2,$ the sequence of SERE
$V_{\varepsilon}(t/\varepsilon) f$ is weakly compact
$\rho - a. s.$ in $D_{B}[0, + \infty)$ with limit points in
$C_{B}[0, + \infty),$ $f \in B_{0}$.

Let us consider the following process in $D_{B}[0, + \infty):$

\begin{eqnarray}\label{e47}
M_{N(t/\varepsilon^{2})}^{\varepsilon} f^{\varepsilon}: = 
V_{N(t/\varepsilon^{2})}^{\varepsilon} f^{\varepsilon} -
f^{\varepsilon} - \sum_{k=0}^{N(t/\varepsilon^{2}) - 1} E_{\rho}
[V_{k+1}^{\varepsilon} f_{k+1}^{\varepsilon} -
V_{k}^{\varepsilon} f_{k}^{\varepsilon}/{\cal F}_{k}],
\end{eqnarray}
where $f^{\varepsilon}: = f + \varepsilon f_{1}
(x(t/\varepsilon^{2})) + \varepsilon^{2} f_{2}
(x(t/\varepsilon^{2}))$, and functions $f_{1}$ and $f_{2}$
are defined from the following equations:

\begin{eqnarray}\label{e48}
(P-I) f_{1}(x) &=& - [m\Gamma(x) + PD_{1}(x, \cdot)] f,\nonumber\\
(P-I) f_{2}(x) &=& [\hat{L} - L(x)] f,\nonumber\\
\hat{L}: &=& \int_{X}\rho(dx)L(x),
\end{eqnarray}

\begin{eqnarray}\label{e272}
L(x) &:=& (m\Gamma(x) + D_{1}(x))
(R_{0} - I)(m\Gamma(x) + D_{1}(x)) +\nonumber\\ 
&\mbox{}& + m_{2}
\Gamma^{2}(x)/2 + mD_{1}(x)
\Gamma(x) + D_{2}(x),\nonumber
\end{eqnarray}

$R_{0}$ is a potential operator of $(x_{n};\quad n \geq 0)$ (see \cite{KS1,KS2,S00}).

The balance condition (\ref{e46}) and condition $\int_X\rho(dx)(\hat{L}-L(x))=0$
give the solvability of the equations in (\ref{e48}).

The process $M_{\nu(t/\varepsilon^{2})}^{\varepsilon} 
f^{\varepsilon}$ is an ${\cal F}_{t}^{\varepsilon}$--martingale
with respect to the $\sigma$--algebra\\
${\cal F}_{t}^{\varepsilon}: = \sigma
\{x(s/\varepsilon^{2}); 0 \leq s \leq t\}$.

This martingale has the asymptotic representation:

\begin{eqnarray}\label{e49}
M_{N(t/\varepsilon^{2})}^{\varepsilon} f^{\varepsilon} =
V_{N(t/\varepsilon^{2})}^{\varepsilon} f - f - \varepsilon^{2}
\sum_{k=0}^{N(t/\varepsilon^{2}) - 1}\hat{L}V_{k}^{\varepsilon} f -
0_{f}(\varepsilon),
\end{eqnarray}
where $\hat{L}$ is defined in (\ref{e48}) and 

$$
\|0_{f}(\varepsilon)\|/\varepsilon \to const\quad \varepsilon \to 0,
\quad \forall f \in B_{0}.
$$


The families $l(M_{N(t/\varepsilon^{2})}^{\varepsilon} 
f^{\varepsilon})$ and $l(\sum_{k=0}^{N(t/\varepsilon^{2})}
E_{\rho}[(V_{k+1}^{\varepsilon} f_{k+1}^{\varepsilon} - 
V_{k}^{\varepsilon} f_{k}^{\varepsilon})/{\cal F}_{k}])$
are weakly compact for all $l \in B_{0}^{*}, \quad f \in B_{0}$.

Set $V^{0}(t)$ for the limit process for $V_{\varepsilon}
(t/\varepsilon)$ as $\quad \varepsilon \to 0$.

Similar to (\ref{e40}) we obtain that the limits for $V^{\varepsilon}
(t)$ and $V_{N(t/\varepsilon^{2})}^{\varepsilon}$
are the some, namely, $V^{0}(t),$  because 

$$
[V^{\varepsilon} (t) - V_{N(t/\varepsilon)}^{\varepsilon^2}] = 
[\Gamma_{x(t/\varepsilon^2)} (t - \varepsilon\tau_{N(t/\varepsilon^2)})
- I] \cdot V_{N(t/\varepsilon^2)}^{\varepsilon}.
$$

Taking into account the LLN for Hawkes process, $N(t)/t\to_{t\to+\infty}\lambda/(1-\hat\mu),$ see e.g., \cite{BDHM,BM,FLM,Z}, the sum $\varepsilon^{2} \sum_{k=0}^{N(t/\varepsilon^{2})}
\hat{L} V_{k}^{\varepsilon} f$ converges strongly
$as \quad \varepsilon \to 0$ to the integral
$\hat\lambda\int_{0}^{t}\hat{L} V^{0}(s)fds,$ where $\hat\lambda$ is defined in (\ref{e41}).

Set $M^{0}(t)f$ be a limit martingale for
$M_{\nu(t/\varepsilon^{2})}^{\varepsilon} f^{\varepsilon}\quad
as \quad \varepsilon \to 0$.

Then, from (32)--(35) and previous reasonings (51)-(53), we have
as $\varepsilon \to 0:$

\begin{eqnarray}\label{e50}
M^{0}(t)f = V^{0}(t)f - f - \hat\lambda\cdot \int_{0}^{t}
\hat{L} V^{0}(s)fds.
\end{eqnarray}

The quadratic variation of the martingale $M^{0}(t)f$ has the
form:

\begin{eqnarray}\label{e51}
<l(M^{0}(t)f)> = \int_{0}^{t}\int_{X}l^{2}
(\sigma\Gamma(x)V^{0}(s)f)\rho(dx)ds,
\end{eqnarray} 
where

$$
\sigma^{2}: = \hat\lambda[m_{2} - m^{2}],
$$
where $m_2, m$ are defined in E), (\ref{e29}), and $\hat\lambda$ in (\ref{e41}).

The solution of martingale problem for $M^{0}(t)$ (namely, to find
the representation of $M^{0}(t)$ with quadratic variation (\ref{e51})) is
expressed by the integral over Wiener orthogonal martingale measure
$W(dx, ds)$ with quadratic variation $\rho(dx)\cdot ds$:

\begin{eqnarray}\label{e52}
M^{0}(t)f = \sigma\int_{0}^{t}\int_{x}\Gamma(x)V^{0}(s)f
W(dx, ds).
\end{eqnarray}

In such a way, the limit process $V^{0}(t)$ satisfies the following
equation (see (\ref{e50}) and (\ref{e52})):

\begin{eqnarray} \label{e53}
V^{0}(t)f = f + \hat\lambda\cdot \int_{0}^{t}\hat{L}\cdot V^{0}(s)
fds + \sigma\int_{0}^{t}\int_{X}\Gamma(x)V^{0}(s)f W(dx, ds).
\end{eqnarray}

If the operator $\hat{L}$ generates the semigroup $U(t)$ then
the process $V^{0}(t)f$ in (\ref{e53}) satisfied equation:

\begin{eqnarray}\label{e54}
V^{0}(t)f = U(t)f + \sigma\int_{0}^{t}\int_{x}U(t-s)
\Gamma(x)V^{0}(s)f W(dx, ds).
\end{eqnarray}

The uniqueness of the limit evolution $V^{0}(t)f$ in 
diffusion approximation scheme follows from the uniqueness 
of the solution of martingale problem for $V^{0}(t)f$ (see
((54)-(56)). The latter is proved by dual SERE in series scheme 
by constructing the limit equation in diffusion approximation 
and by using a dual identify. For more details, see \cite{KS1}. {\it Q.E.D}.

{\bf Remark.} In a similar way, we can consider diffusion approximation of self-exciting random evolutions in 
reducible phase space, as we did it in \cite{KS1,KS2,S00} for semi-Markov RE.







\section{Discussion}
The main difference between SERE and, e.g., semi-Markov REs  (SMRE) is that it involves new features for the  REs: self-exciting and clustering ones, which is not the case for the SMRE.
This was incorporated into the Swish process $x(t):=x_{N(t)},$ where $N(t)$ is a Hawkes process, and into the parameter $\hat\lambda:=\frac{\lambda}{1-\hat\mu},$ $\hat\mu:=\int_0^{+\infty}\mu(t)dt=\alpha/\beta,$ see, e.g., sec. 6.1, for the averaged self-exciting summation on a Markov chain. This parameter is presented in all results for the SERE, including averaging and diffusion approximation, and their applications. These features show originality and novelty of the paper. The results for SERE obtained in the paper generalize many results in the literature devoted to the Hawkes-based models in finance and insurance, such as, e.g.,  risk model based on general compound Hawkes process \cite{S2} (self-exciting impulse traffic/transport process (8) in this paper) and general compound Hawkes processes in limit order books \cite{S8} (self-exciting summation on a Markov chain in this paper). In applications, such as high-frequency and algorithmic trading (HFT) and in limit order books modelling (when transactions happen every milliseconds), one usually considers the long scale $tn$ (minutes, hours, etc.) instead of $t$ (milliseconds) to have time for solutions many problems like liquidation, acquisitions, and market making, to name a few. In this case, the limit theorems such as law of large numbers (averaging) and functional central limit theorem (diffusion approximation) for the case when $n\to+\infty$ and $t$ is fixed, are very useful. In this paper, we used the scaling $t/\epsilon$ to get our results when $\epsilon\to 0.$ Of course, it is similar to the case $tn,$ when $n\to+\infty$ (see, e.g., \cite{S8}). Therefore, many applications of SERE already exist in finance  \cite{S8} and risk theory \cite{S2}. However, many more applications of SERE to the examples of SERE in sec. 3  could be considered space permitting for the current paper.  Finally, we note that some new results such as normal deviations and rates of convergence for SERE and their applications could also be formulated space permitting for the paper. Therefore, the author leaves  preparing those results and applications to the future research work. 

\section{Conclusion }
In this paper we have developed  a new class of so-called self-exciting random evolutions (SEREs) and their applications. We have first introduced a new random process $x(t)$ such that it is based on  a superposition of a Markov chain $x_n$ and a Hawkes process $N(t),$ i.e., $x(t):=x_{N(t)}.$ We call this process self-walking imbedded semi-Hawkes process (Swish Process or SwishP). Then the self-exciting REs (SEREs) have been constructed in similar way as, e.g., semi-Markov REs, but instead of semi-Markov process $x(t)$ we have SwishP.  We gave classifications and examples of self-exciting REs (SEREs). Then we considered two limit theorems for SEREs such as  averaging  (Theorem 1) and diffusion approximation (Theorem 2). Applications of SEREs have been devoted to the so-called self-exciting traffic/transport process and self-exciting summation on a Markov chain, which are examples of continuous and discrete SERE. From these processes we can constructed many other self-exciting processes, e.g., such as impulse traffic/transport process, self-exciting risk process, etc. We also presented averaged and diffusion approximation for some self-exciting processes.

\section*{Acknowledgment} The author thanks NSERC for continuing support.


\end{document}